\documentclass[12pt]{amsart}

\usepackage{amsfonts}

\usepackage{float}

\usepackage{epsfig}

\usepackage{amscd}

\usepackage{latexsym}

\newcommand{\square}{\Box}
\newtheorem{theorem}{Theorem}
\newtheorem{lemma}[theorem]{Lemma}
\newtheorem{proposition}[theorem]{Proposition}
\newtheorem{corollary}[theorem]{Corollary}
\newtheorem{definition}[theorem]{Definition}
\newtheorem{example}[theorem]{Example}
\newtheorem{conjecture}[theorem]{Conjecture}
\newtheorem{remark}[theorem]{Remark}

\begin{document}

\author{J\o rgen Ellegaard Andersen}
\address{Dept. of Math., Univ. of Aarhus, DK-8000 Aarhus C, DENMARK}
\curraddr{MSRI, CA-94720, USA}

\email{andersen@msri.org}

\author{Josef Mattes}
\address{Dept. of Math., UC Davis, CA-95616,
USA, http://math.ucdavis.edu/\symbol{126}mattes}

\author{Nicolai Reshetikhin}
\address{Dept. of Math., UC Berkeley, CA-94720, USA}

\title{Quantization of the Algebra of Chord Diagrams}

\thanks{Joint MSRI and University of Aarhus Preprint.}
\thanks{Research at MSRI is supported in part by NSF grant DMS-9022140.}
\thanks{J.E.A. is supported in part by NSF grant DMS-93-09653 and by the
Danish Research Council.}

\begin{abstract}
In this paper we define an algebra structure on the vector space $L(\Sigma )$
generated by links in the manifold $\Sigma \times [0,1]$ where $\Sigma $ is
an oriented surface. This algebra has a filtration and the associated graded
algebra $L_{Gr}(\Sigma )$ is naturally a Poisson algebra. There is a Poisson
algebra homomorphism from the algebra of chord diagrams $ch\left( \Sigma
\right) $ on $\Sigma $ to $L_{Gr}(\Sigma )$.

We show that multiplication in $L\left( \Sigma \right) $ provides a
geometric way to define a deformation quantization of the algebra of chord
diagrams, provided there is a universal Vassiliev invariant for links in $%
\Sigma \times [0,1]$. The quantization descends to a quantization of the
moduli space of flat connections on $\Sigma $ and it is universal with
respect to group homomorphisms. If $\Sigma $ is a compact with free
fundamental group we construct a universal Vassiliev invariant.
\end{abstract}

\maketitle

\tableofcontents

\section{Introduction}

The study of finite type (or Vassiliev) invariants of knots and links \cite
{t8}, \cite{vassiliev92} leads to the notion of chord diagrams \cite
{Bar-Natan92}. All the known quantum invariants of links are determined by
finite type invariants \cite{t1}.

In \cite{Mattes96a} we extended the notion of a chord diagram to arbitrary
oriented Riemann surfaces $\Sigma $ and showed that these chord diagrams are
an effective tool to clarify the study of the algebra of functions on the
moduli space of flat connections on $\Sigma $ as started by Goldman in \cite
{m1}. The algebra of chord diagrams $ch\left( \Sigma \right) $ is a graded
Poisson algebra, its component of degree zero is exactly the algebra of
loops considered in \cite{m1}. There is a natural Poisson algebra
homomorphism from the Poisson algebra of (coloured) chord diagrams to $%
\mathcal{F}\left( \mathcal{M}_{\Sigma }^{G}\right) $, the functions on the
moduli space of flat $G$-connections, which for many simple groups $G$ is
known to be surjective \cite{Mattes96a},\cite{Mattes96b}. Therefore $%
ch\left( \Sigma \right) $ can be viewed as a universal moduli space.

The present paper is the first in a series devoted to the study of the
quantum versions of the structures used in \cite{Mattes96a}. We show that
there is a universal invariant of links in $\Sigma \times \left[ 0,1\right] $
with values in chord diagrams, where $\Sigma $ is a punctured surface. Using
this we construct a quantization of the algebra of chord diagrams on $\Sigma 
$ using the (noncommutative) multiplication of links in $\Sigma \times
\left[ 0,1\right] $. We show that (in contrast to the ''obvious'' linear
quantization of $ch\left( \Sigma \right) $) our quantization descends to a
quantization of the moduli space $\mathcal{M}_{\Sigma }^{G}$ of flat $G$%
-connections on $\Sigma $.

Other approaches to the quantization of the moduli space of flat connections
include skein modules as in \cite{Turaev91},\cite{Ozawa} and the related
quantization of the $SL_{2}$-moduli space of \cite{Bullock96}, geometric
quantization \cite{atiyah90}, \cite{Axelrod91}, \cite{Faltings93}, \cite
{Hitchin90} and the combinatorial approach of \cite{Alekseev95a}, \cite
{Alekseev95c}, \cite{Alekseev95}. We will comment on relations to these
constructions in a forthcoming publication.

Chord diagrams also arose in Yang-Mills
theory on Minkowski space \cite{Rajeev95a} and in quantum gravity \cite
{Baez1994b}.

We have been greatly helped by discussions with Dylan Thurston and Pol
Vanhaecke as well as by Thang Le, Pierre Vogel, Akira Yoshioka and by Dror Bar-Natan's
Mathematica package for calculating with chord diagrams, available at \newline
http://www.ma.huji.ac.il/\symbol{126}drorbn, cf. \cite{Bar-Natan93}.

J.E.A. would like to thank the Department of Mathematics, University
of California, Berkeley and the Mathematical Sciences Research Institut for their
hospitality during the completion of this work.

\section{Algebra of links in $\Sigma \times \left[ 0,1\right] $}

\subsection{Filtered and graded rings}

By a filtered ring we mean a ring $F$ with a descending filtration $%
F=F_{0}\supseteq F_{1}\supseteq \dots $ indexed by the natural
 numbers,
such that the multiplication respects the filtration $F_{m}F_{n}\subseteq
F_{m+n}$. We write $F_{\infty }=\bigcap_{n\geq 0}F_{n}$ and $F^{\prime
}=F/F_{\infty }$. A map of filtered rings $f:F\rightarrow G$ has to respect
the filtrations: $f\left( F_{n}\right) \subseteq G_{n}$.

Each filtered ring defines a graded ring $F_{Gr}=\bigoplus_{n\geq
0}F^{\left( n\right) }$ with graded components $F^{\left( n\right)
}=F_{n}/F_{n+1}$. A map of filtered rings $f$ induces a map $f_{Gr}$ of the
associated graded rings. The associative multiplication on $F$ induces an
associative multiplication on $F_{Gr}$: For $x\in F_{n},y\in F_{m}$ define
the product of their classes $[x]\in F^{(n)},[y]\in F^{(m)}$ as the class of
the product: 

\begin{equation}
\lbrack x].[y]:=[xy]\in F^{(n+m)}  \label{quant3}
\end{equation}
If $f$ is a ring homomorphism, so is $f_{Gr}$. Conversely, given a graded ring
$G=\bigoplus_{n\geq 0}G^{\left( n\right) }$ defines a filtered ring $G^{F}$ with
filtered components $ G_{m}^{F}=\bigoplus_{n\geq m}G^{\left( n\right) }$ and
another filtered ring $\overline{G}$ with filtered components
$\overline{G}_{m}=\prod_{n\geq m}G^{\left( n\right) }$. This is the topological
completion of $G^{F}$ in the topology whose basis of open sets is given by
translates of the filtered components. Clearly $F^{\prime }\subseteq
\overline{F_{Gr}}$.

\begin{lemma}

\label{a211096}If $F=F_{0}\supseteq F_{1}\supseteq F_{2}\supseteq \dots $ is
a filtered ring such that the associated graded ring $F_{Gr}=%
\bigoplus_{n=0}^{\infty }F_{n}/F_{n+1}$ is commutative, then the bracket 
\begin{equation}
\{[x],[y]\}:=[xy-yx]\in F^{(n+m+1)}  \label{quant2}
\end{equation}
for $x\in F_{n},y\in F_{m}$ defines a Poisson bracket on $%
F_{Gr}$.

\end{lemma}

\subsection{Multiplication of links}

If $M_{1}$ and $M_{2}$ are two 3-manifolds, $X\hookrightarrow \partial M_{1}$
and $X^{-}\hookrightarrow \partial M_{2}$, then we denote by $M_{1}\cup
_{X}M_{2}$ the result of gluing $M_{1}$ and $M_{2}$ along $X$. Denote by $%
L(\Sigma )$ the ring spanned by (framed) links in $\Sigma \times [0,1]$ (by
a link we mean an isotopy class of smooth imbeddings $(S^{1})^{\cup
k}\hookrightarrow \Sigma \times [0,1]$), with multiplication in $L(\Sigma )$
defined by the following composition of maps:

\[
L_1\otimes L_2\mapsto L_1\times L_2\mapsto L_1\cup L_2 
\]

where $L_1\times L_2\subset (\Sigma \times [0,1/2])\times (\Sigma \times
[1/2,1])$ and $L_1\cup L_2\subset (\Sigma \times [0,1/2])\cup _{\Sigma
\times \left\{ 1/2\right\} }(\Sigma \times [1/2,1])=\Sigma \times [0,1]$ .

\begin{proposition}

This multiplication determines on $L(\Sigma )$ the structure of an
associative (in general noncommutative) ring with the empty link being the
unit element.

\end{proposition}

%TCIMACRO{\TeXButton{Proof}{\proof}}

%BeginExpansion

{\bf Proof:} Obvious. 
%TCIMACRO{\TeXButton{End Proof}{\endproof}}
%BeginExpansion
$\square$%

%EndExpansion

\subsection{Filtration\label{quant0}}

Let $L\subset \Sigma \times [0,1]$ be a link and $D_{L}\subset \Sigma $ some
link diagram of $L$, so that $D_{L}$ is (an isotopy class of) a regular
projection of $L$ to $\Sigma $. As usual we distinguish vertices of two
types: For each vertex $v\in D_{L}$ we introduce an oriented crossing number 
$\epsilon (v)$ as the sign of the vertex

$+$ is 
%TCIMACRO{
%\TeXButton{quant1fig3.eps}{\epsfig{file=q1fig3.eps,height=3cm,width=3cm,angle=-90}
%}}
%BeginExpansion
\epsfig{file=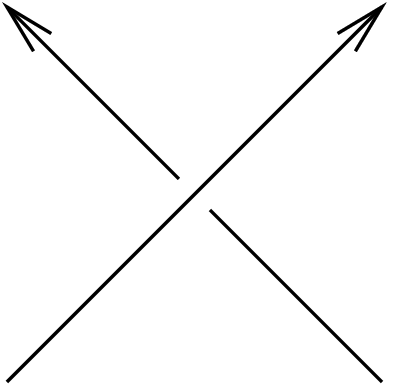,height=1cm,width=1cm,angle=-90}
%
%EndExpansion
, $-$ is 
%TCIMACRO{
%\TeXButton{quant1fig4.eps}{\epsfig{file=q1fig4.eps,angle=-90,height=2cm,width=2cm}
%}}
%BeginExpansion
\epsfig{file=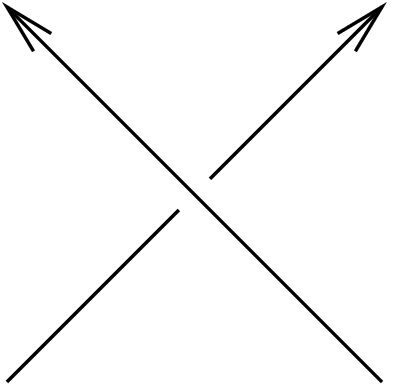,height=1cm,width=1cm,angle=-90}
%

%EndExpansion

If the diagram $D_{L}$ has vertices $v_{1},\dots ,v_{n}$ with corresponding
oriented crossing numbers $\epsilon _{1},\dots ,\epsilon _{n}$ we also
denote it by $D_{L}^{\epsilon _{1},\dots ,\epsilon _{n}}$ when we wish to
emphasize the types of the crossings. We may regard $L$ as an equivalence
class $\left[ D_{L}\right] $ of diagrams that are related by Reidemeister
moves.

Introduce the following operation $\nabla $: Choose a set of crossings $%
v_{i_{1}},\dots ,v_{i_{m}}$ of $D_{L}$ and set 

\[
\nabla _{v_{i_{1}},\dots ,v_{i_{m}}}D_{L}:=\sum_{\epsilon _{i_{1}},\dots
,\epsilon _{i_{m}}=\pm 1}\epsilon _{i_{1}}\dots \epsilon
_{i_{m}}[D_{L}^{\epsilon _{1},\dots ,\epsilon _{n}}] 
\]

This maps the link diagram $D_{L}$ to a linear combination of links whose
regular projections are obtained from $D_{L}$ by switching crossings. If $%
v_{1},\dots ,v_{m}$ is a list of all crossing of $D_{L}$ then we denote $%
\nabla _{v_{i_{1}},\dots ,v_{i_{m}}}D_{L}$ simply by $\nabla D_{L}$.

Now let $L_{m}\left( \Sigma \right) \subseteq L\left( \Sigma \right) $ be
the span of all elements of the form $\nabla _{v_{i_{1}},\dots
,v_{i_{m}}}D_{L}$ where $D_{L}$ runs over all possible link diagrams.

\begin{proposition}

\label{a260996}

\begin{enumerate}

\item  The filtration $L(\Sigma )\supset L_{1}(\Sigma )\supset L_{2}(\Sigma
)\dots $ is compatible with the algebra structure.

\item  If $x\in L_{n}$ and $x^{\prime }\in L_{m}$ then $x.x^{\prime
}-x^{\prime }.x\in L_{m+n+1}$.

\end{enumerate}

\end{proposition}

%TCIMACRO{\TeXButton{Proof}{\proof}}

%BeginExpansion

{\bf Proof :}%

%EndExpansion

\begin{enumerate}

\item  We have to prove that $L_{n}(\Sigma ).L_{m}(\Sigma )\subset
L_{n+m}(\Sigma )$. Indeed, on generators, if $x=\nabla _{v_{1},\dots
,v_{n}}D_{L},x^{\prime }=\nabla _{w_{1},\dots ,w_{m}}D_{L^{\prime }}$ are in
general position, then $x.x^{\prime }=\nabla _{v_{1},\dots
,v_{n},w_{1},\dots ,w_{m}}D_{L.L^{\prime }}$.

\item  We know that $x.x^{\prime }-x^{\prime }.x=\nabla _{v_{1},\dots
,v_{n},w_{1},\dots ,w_{m}}D_{L.L^{\prime }}-\nabla _{v_{1},\dots
,v_{n},w_{1},\dots ,w_{m}}D_{L^{\prime }.L}$ and one can obtain $%
D_{L^{\prime }.L}$ from $D_{L.L^{\prime }}$ by a sequence of crossing
changes away from $v_{1},\dots ,v_{n},w_{1},\dots ,w_{m}$. This clearly
implies the assertion.

\end{enumerate}
%TCIMACRO{\TeXButton{End Proof}{\endproof}}
%BeginExpansion
$\square$%

%EndExpansion

It follows from proposition \ref{a260996}(2) that the multiplication on $%
L_{Gr}$ is commutative. By continuity it extends to $\overline{L_{Gr}\left(
\Sigma \right) }$.

\subsection{Finite type invariants}

Recall that every link invariant extends in the obvious way to an element of 
$L\left( \Sigma \right) ^{*}$, we identify the invariant with this
extension. A link invariant $v$ is called a finite type invariant (or
Vassiliev invariant) if $v|_{L_{n}\left( \Sigma \right) }=0$ for some $n$.
The degree of $v$ is the smallest such $n$ minus $1$.

Set $L_{\infty }\left( \Sigma \right) :=\bigcap_{n\in \Bbb{N}}L_{n}\left(
\Sigma \right) $. We will be interested mostly in the algebra $L^{\prime
}\left( \Sigma \right) :=L\left( \Sigma \right) /L_{\infty }\left( \Sigma
\right) $. The question whether $L\left( \Sigma \right) =L^{\prime }\left(
\Sigma \right) $ is difficult, in fact it is equivalent to the main open
problem in the theory of finite type invariants: whether finite type
invariants are sufficient to determine link types.

\begin{proposition}

Finite type invariants distinguish links in $\Sigma \times \left[ 0,1\right] 
$ iff $L_{\infty }\left( \Sigma \right) =\emptyset $.

\end{proposition}

%TCIMACRO{\TeXButton{Proof}{\proof}}

%BeginExpansion

{\bf Proof:} 
If $L_{\infty }\left( \Sigma \right) =\emptyset $, then for any pair of
distinct links $L_{1},L_{2}$ there is some $n$ such that $\left[
L_{1}-L_{2}\right] \neq 0\in L\left( \Sigma \right) /L_{n}\left( \Sigma
\right) $. Therefore there is an element $v$ of $\left( L\left( \Sigma
\right) /L_{n}\left( \Sigma \right) \right) ^{*}$, i.e. an invariant of type 
$\leq n$, such that $v\left( L_{1}\right) -v\left( L_{2}\right) \neq 0$.

Conversely, assume $\sum_{i=0}^{n}\lambda _{i}L_{i}\in L_{\infty }\left(
\Sigma \right) $. By induction on $n$ one constructs a finite type invariant 
$v$ such that $v\left( L_{i}\right) \neq v\left( L_{j}\right) $ for $i\neq j$%
. Using that any polynomial in $v$ is again a finite type invariant and that 
$\det \left( v^{k}\left( L_{l}\right) \right) _{k,l}\neq 0$, we can clearly
find a finite type invariant $w=p\left( v\right) $ with $\sum_{i=0}^{n}%
\lambda _{i}w\left( L_{i}\right) \neq 0$. Hence $\sum_{i=0}^{n}\lambda
_{i}L_{i}\neq 0$ in $L\left( \Sigma \right) /L_{\deg w}\left( \Sigma \right) 
$ which implies $\sum_{i=0}^{n}\lambda _{i}L_{i}\notin L_{\deg w}\left(
\Sigma \right) \supseteq L_{\infty }\left( \Sigma \right) $.%
%TCIMACRO{\TeXButton{End Proof}{\endproof}}
%BeginExpansion
$\square$%

%EndExpansion

\begin{remark}

We could not find this result in the literature, even for $\Sigma =S^{2}$.

\end{remark}

\section{Chord diagrams and universal finite type invariants}

\subsection{Chord diagrams}

Now let us give an explicit description of the vector spaces $L_{n}(\Sigma )$%
. Recall the following definitions from \cite{Mattes96a}:

\begin{definition}

A \emph{chord diagram} is a graph consisting of disjoint oriented circles $%
S_{i},i\in \{1,...,n\}$ and disjoint arcs $C_{j},j\in \{1,...,m\}$ such that:

\begin{enumerate}

\item  the endpoints of the arcs are distinct

\item  $\cup _{j}\partial C_{j}=\left( \cup _{i}S_{i}\right) \cap \left(
\cup _{j}C_{j}\right) $

\end{enumerate}

The arcs are called \emph{chords}, the circles $S_{i}$ are called the \emph{%
core components} of the diagram.

\end{definition}

\begin{definition}

Given a closed oriented surface $\Sigma $, a \emph{geometrical chord diagram}
on $\Sigma $ is a smooth map from a chord diagram $D$ to $\Sigma $, mapping
the chords to points. A \emph{chord diagram on}\textbf{\ $\Sigma $} is a
class of geometric chord diagrams modulo homotopy.

\end{definition}

Note that diffeomorphism classes of chord diagrams correspond to chord
diagrams on $S^{2}$.

\begin{definition}

By a \emph{generic chord diagram} (on $\Sigma $) we will mean a geometrical
chord diagram on $\Sigma $ such that all circles are immersed, and with all
double points transverse.

\end{definition}

Clearly every chord diagram on $\Sigma $ contains generic chord diagrams. A
generic chord diagram on $\Sigma $ with images of chords at points $%
v_{1},\dots ,v_{m}$ will be denoted by $D\left( v_{1},\dots ,v_{m}\right) $.

Consider the complex vector space $V_{\Sigma }$ with the basis given by the
set of chord diagrams on $\Sigma $ and the subspace $W$ generated by the
4T-relations \cite{Bar-Natan92}.

\begin{definition}

The algebra $ch\left( \Sigma \right) :=V_{\Sigma }/W$ is called the \emph{algebra
of chord diagrams on }$\Sigma $.

\end{definition}

It has a natural ring structure with multiplication given by union of chord
diagrams, with unit the empty diagram. It is isomorphic to the polynomial
ring on the space of connected chord diagrams.

These rings are graded by the number of chords 
\[
ch\left( \Sigma \right) =\bigoplus_{n\geq 0}ch^{(n)}\left( \Sigma \right) 
\]
and we have an associated filtered space with filtered components $%
ch_{m}\left( \Sigma \right) :=\bigoplus_{n\geq m}ch^{(n)}\left(
\Sigma \right) $ and completion $\overline{ch\left( \Sigma \right) }%
=\prod_{n\geq 0}ch^{(n)}\left( \Sigma \right) $.

To any element $D\in ch^{\left( n\right) }\left( \Sigma \right) $ we can
associate an element $\lambda \left( D\right) \in L\left( \Sigma
\right) /L_{n+1}\left( \Sigma \right) $ by setting 
\[
\lambda \left( D\right) :=\nabla _{v_{i_{1}},\dots ,v_{i_{n}}}D_{L}\text{ mod }%
L_{n+1}\left( \Sigma \right) 
\]
for any link $L$ that projects to the diagram $D$, where $\left\{
v_{i_{1}},\dots ,v_{i_{n}}\right\} $ is the set of chords of $D$. This defines
a graded linear map $\lambda : ch(\Sigma) \rightarrow L_{Gr}(\Sigma)$.

Given a Lie group $G$ with invariant inner product one can generalize all of
the above in a trivial way to define an algebra $ch^{G}\left( \Sigma \right) 
$ of coloured chord diagrams where each core component is coloured by a
finite dimensional representation $V$ of $G$.

We will also need the notion of \emph{chord tangles} which is just that of a
chord diagram on $\Bbb{R}\times \left[ 0,1\right] $ except that we allow the 
$S_{i}$ to be neatly imbedded intervals rather than circles (so that the
image of $\cup S_{i}$ forms a tangle rather than a link). These chord
tangles will be related to tangles as chord diagrams are to links. Chord
tangles with suitable numbers and orientations of endpoints can be composed
in the obvious way.

\subsection{The Poisson structure}

Recall from \cite{Mattes96a} that $ch\left( \Sigma \right) $ has a natural
Poisson structure given as follows: Assume $D_{1}\cup D_{2}$ is a generic
chord diagram. For $p\in D_{1}\cap D_{2}$ we define the oriented
intersection number by $\epsilon _{12}(p):=$ $\left\{ 
\begin{array}{cc}
+ & \text{for 
%TCIMACRO{
%\TeXButton{picture4}{\begin{picture}(30,20)\put(5,0){\vector(1,1){20}}\put(25,0){\vector(-1,1){20}}\put(0,10){\thinlines 1}\put(30,10){\thinlines 2}\put(15,0){\thinlines p}\end{picture}} }
%BeginExpansion
\begin{picture}(30,20)\put(5,0){\vector(1,1){20}}\put(25,0){\vector(-1,1){20}}\put(0,10){\thinlines 1}\put(30,10){\thinlines 2}\put(15,0){\thinlines p}\end{picture}%
%EndExpansion
} \\ 
- & \text{for 
%TCIMACRO{
%\TeXButton{picture5}{\begin{picture}(30,20)\put(5,0){\vector(1,1){20}}\put(25,0){\vector(-1,1){20}}\put(0,10){\thinlines 2}\put(30,10){\thinlines 1}\put(15,0){\thinlines p}\end{picture}} }
%BeginExpansion
\begin{picture}(30,20)\put(5,0){\vector(1,1){20}}\put(25,0){\vector(-1,1){20}}\put(0,10){\thinlines 2}\put(30,10){\thinlines 1}\put(15,0){\thinlines p}\end{picture}%
%EndExpansion
}
\end{array}
\right. $ \hspace{0.3cm} where 1 and 2 indicate components of the
corresponding diagrams.

For each for $p\in D_{1}\cap D_{2}$ we define $D_{1}\cup _{p}D_{2}$ to be
the chord diagram on $\Sigma $ given by joining $D_{1}^{-1}(p)\in $ and $%
D_{2}^{-1}(p)$ by a chord. Under the above assumptions, for chord diagrams $%
D_{1},D_{2}$ we define their Poisson bracket to be

\begin{equation}
\left\{ \lbrack D_{1}],[D_{2}]\right\} :=\sum_{p\in D_{1}\cap D_{2}}\epsilon
_{12}(p)[D_{1}\cup _{p}D_{2}]  \label{b270195}
\end{equation}

\begin{proposition}

\label{b260996}The map $\lambda :ch\left( \Sigma \right) \rightarrow
L_{Gr}\left( \Sigma \right) $ is a homomorphism of graded Poisson rings.
\end{proposition}

%TCIMACRO{\TeXButton{Proof}{\proof}}

%BeginExpansion

{\bf Proof:} Given two chord diagrams $D_{1},D_{2}$ in generic position with $D_{1}\cap
D_{2}=\{p_{n+m+1},\dots ,p_{n+m+k}\}$ the Poisson structures compare as
follows: 
\begin{eqnarray*}
\left\{ \lambda \left( D_{1}\right) ,\lambda \left( D_{2}\right) \right\}
&=&\left[ \nabla _{v_{1},\dots ,v_{n}}D_{L_{1}},\nabla _{w_{1},\dots
,w_{m}}D_{L_{2}}\right] \text{ mod }L_{n+1+m+1}\left( \Sigma \right) \\
&=&\nabla _{v_{1},\dots ,v_{n},w_{1},\dots ,w_{m}}\left(
D_{L_{1}.L_{2}}-D_{L_{2}.L_{1}}\right) \\
&=&\nabla _{v_{1},\dots ,v_{n},w_{1},\dots ,w_{m}}\left(
\sum_{i=1}^{k}\varepsilon _{i}\nabla _{p_{i}}D_{L_{1}.L_{2}}\right) \text{mod%
}L_{n+m+2}\left( \Sigma \right) \\
&=&\sum_{i=1}^{k}\varepsilon _{i}\nabla _{v_{1},\dots ,v_{n},w_{1},\dots
,w_{m},p_{i}}D_{L_{1}.L_{2}} \\
&=&\lambda \left( \left\{ D_{1},D_{2}\right\} \right)
\end{eqnarray*}
where we used 

\begin{eqnarray*}
D_{L_{1}.L_{2}}-D_{L_{2}.L_{1}} &=&D_{L_{1}.L_{2}}^{\varepsilon _{1},\dots
,\varepsilon _{n+m},\varepsilon _{n+m+1},\dots ,\varepsilon
_{n+m+k}}-D_{L_{1}.L_{2}}^{\varepsilon _{1},\dots ,\varepsilon
_{n+m},-\varepsilon _{n+m+1},\dots ,-\varepsilon _{n+m+k}} \\
&=&\sum_{i=1}^{k}\left( D_{L_{1}.L_{2}}^{\varepsilon _{1},\dots ,\varepsilon
_{n+m+i-1},\varepsilon _{n+m+i},-\varepsilon _{n+m+i+1},\dots ,-\varepsilon
_{n+m+k}}\right. \\
&&\left. -D_{L_{1}.L_{2}}^{\varepsilon _{1},\dots ,\varepsilon
_{n+m+i-1},-\varepsilon _{n+m+i},-\varepsilon _{n+m+i+1},\dots ,-\varepsilon
_{n+m+k}}\right) \\
&=&\sum_{i=1}^{k}\varepsilon _{i}\nabla _{p_{i}}D_{L_{1}.L_{2}}^{\varepsilon
_{1},\dots ,\varepsilon _{n+m+i},-\varepsilon _{n+m+i+1},\dots ,-\varepsilon
_{n+m+k}} \\
&=&\sum_{i=1}^{k}\varepsilon _{i}\nabla _{p_{i}}D_{L_{1}.L_{2}}\text{ mod }%
L_{n+m+2}\left( \Sigma \right)
\end{eqnarray*}
%TCIMACRO{\TeXButton{End Proof}{\endproof}}
%BeginExpansion
$\square$
%EndExpansion

This Poisson structure trivially extends to the case of coloured chord
diagrams. It is closely related to the Poisson structure on the moduli space
of flat connections on $\Sigma $. The following is one of the main results
of \cite{Mattes96a}:

\begin{theorem}

There is a Poisson algebra homomorphism from $ch^{G}\left( \Sigma \right) $
to the Poisson algebra $\mathcal{F}\left( \mathcal{M}^{G}\right) $ of
functions on the moduli space of flat $G$-connections on $\Sigma $. This
homomorphism is universal with respect to Lie group homomorphisms that
preserve the invariant inner product.
\end{theorem}

This Poisson algebra homomorphism known to be surjective for many
interesting groups (\cite{Mattes96a},\cite{Mattes96b}). We will see later
that this allows us to obtain quantizations of the algebra of functions on
moduli space from our quantization of chord diagrams.

\subsection{Vassiliev-Kontsevich invariants}

Recall that for the 2-sphere there exists the notion of a \emph{universal
Vassiliev invariant} (or universal Vassiliev-Kontsevich invariant \cite
{Kontsevich},\cite{Bar-Natan92}). It can be regarded as an invariant of
links in the cylinder $S^{2}\times \left[ 0,1\right] $ that takes values in
the algebra $\overline{ch\left( S^{2}\right) }$. If $L$ is a link in $%
S^{2}\times \left[ 0,1\right] $ then this can be written as a formal sum
over all chord diagrams 
\[
V\left( L\right) =\sum_{D\in ch(S^2)} \langle D,L\rangle D. 
\]
Here $\langle D,L\rangle $ should be certain intersection numbers which
can be written in many different ways (see \cite{Kontsevich,BT,DT,PV}).
Similarly, for any compact oriented surface
$%
\Sigma $ we expect that there is an analogous universal invariant of links in $%
\Sigma \times \left[ 0,1\right] $ with values in $ch(\Sigma )$,
\[
V\left( L\right) =\sum_{D\in ch(\Sigma)} \langle D,L\rangle D,
\]
having the following important property: $V\left( L_{n}\left( \Sigma \right)
\right) \subseteq \bigoplus_{m\geq n}ch^{\left( m\right) }\left( \Sigma
\right) $ and for any chord diagram $D\in ch^{\left( k\right) }\left( \Sigma
\right) $%

\begin{equation}
V\left( \lambda \left( D\right) \right) =D\text{ }\text{ mod }ch^{\left(
k+1\right) }\left( \Sigma \right)  \label{a230896}
\end{equation}

The existence of a universal Vassiliev invariant is also known as the \emph{%
fundamental theorem for Vassiliev invariants} \cite{Bar-Natan96}.
The approach that we will use follows the construction of
\cite{Bar-Natan93} closely: Given a link $L$, the link projection is decomposed
into a composition of elementary tangles, then to each tangle one associates a
series of string chord diagrams (in the terminology of \cite{Mattes96a}) and
the invariant $V\left( L\right) $ is the composition of the string chord
diagrams, i.e. a series of chord diagrams in $\overline{ch\left( \Sigma
\right) }$.

\begin{figure}

\begin{center}

\epsfig{file=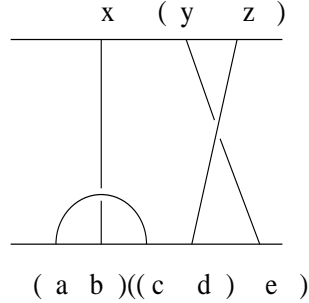,height=4cm,width=4cm,angle=-90}

\caption{A possible non-associative tangle with upper endpoints $x(yz)$ and lower
endpoints $(ab)((cd)e)$\label{quant1fig5}}

\end{center}

\end{figure}

We recall the construction from \cite[definition 2.5]{Bar-Natan93}. A \emph{%
non-associative tangle} is a tangle where the upper and lower endpoints are
parenthesized, see figure \ref{quant1fig5}.

\begin{figure}

\begin{center}

\epsfig{file=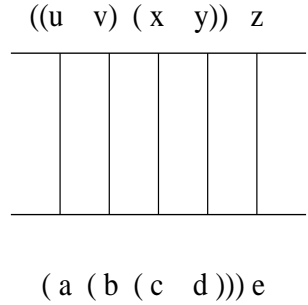,height=4cm,width=4cm,angle=-90}

\caption{A possible associativity morphism $(ab)((cd)e)\rightarrow
((uv)(xy))z$\label{quant1fig6}}

\end{center}

\end{figure}

\begin{figure}
\begin{center}
\epsfig{file=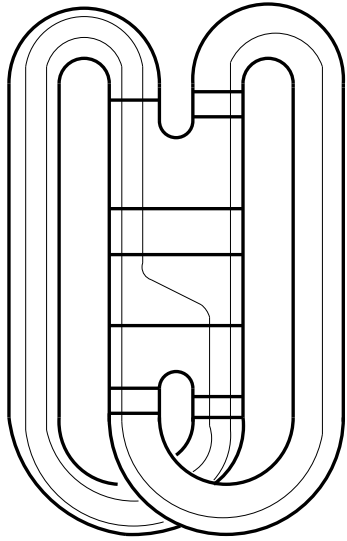,angle=-90}
\caption{A possible decomposition of a diagram on the punctured torus
\label{quant1fig1}}
\end{center}
\end{figure}

Clearly, nonassociative tangles with suitable bracketings and orientations
of endpoints can be composed. Every non-associative tangle can be decompose
into elementary ones: crossings (denoted G2 in Bar-Natan's paper), local
maxima and minima (G3), and trivial tangles $\text{id}$, where only one
parenthesis is moved in the form $\left( A\left( BC\right) \right)
\longleftrightarrow \left( \left( AB\right) C\right) $ (\emph{associativity
morphisms}, G1), see for example figure \ref{quant1fig6}. Associativity
morphisms are often drawn so that the location of parenthesis is indicated
by the distance between strands in the tangle. An example of this can be
seen in figure \ref{quant1fig1}.

It follows from Drinfeld's work \cite{Drinfeld}, that one can associate formal
series of chord tangles $\Phi ,R,e,i$ to the elementary nonassociative tangles
G1,G2,G3 such that the following holds: If $U\left( T_{i}\right) $ denotes the
series corresponding to the elementary tangle $T_{i}$ and a given nonassociative
tangle $T$ is decomposed into elementary ones $T=\prod_{i=1}^{n}T_{i}$ then

\begin{theorem}[Drinfeld,Kontsevich,Bar-Natan]

\label{b161096}The series of chord tangles $\prod_{i=1}^{n}U\left(
T_{i}\right) $ is an invariant of $T$.

\end{theorem}

This implies the fundamental theorem of Vassiliev invariants for $S^2\times
[0,1]$ (by letting $T$ have no endpoints, hence no parenthesis). Analogously we
can show

\begin{theorem}

\label{a220796}If $\Sigma $ is a compact Riemann surface with nonempty
boundary, then there exists a universal Vassiliev-Kontsevich invariant.

\end{theorem}

{\bf Proof:} Without loss of generality we can assume $\Sigma $ is connected and not
equal to $\Bbb{D}^{2}$.  Given a link in $\Sigma\times [0,1]$, we can cut 
$\Sigma$ along a
finite number of neatly embedded intervals $I_j$, to obtain a
decomposition $\Sigma = \cup_i \Sigma^S_i \cup \cup_k \Sigma^H_k$ where
the $\Sigma^S_i$'s are squares and the
$\Sigma^H_k$ are hexagons, with the following properties: We can choose a
bracketing of all the ordered sets $L\cap I_j$ such that

\begin{figure}
\begin{center}
\epsfig{file=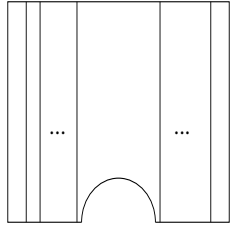,height=3cm,width=3cm,angle=-90}
\caption{The hexagon non-associative tangle G5
\label{G5}}
\end{center}
\end{figure}

\begin{description}
\item[1.] For each square $\Sigma^S_i$, $L^S_i=\Sigma^S_i\cap L \subset
\Sigma^S_i\times [0,1]$ with the bracketing of the top and bottom of $L^S_i$ is an
elementary non-associative tangle of type G1, G2, G3 or a trivial
non-associative tangle G4 as in figure \ref{G4}, where the bracketing is the same
at the top and at the bottom.
\begin{figure}
\begin{center}
\epsfig{file=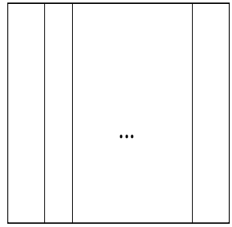,height=3cm,width=3cm,angle=-90}
\caption{The trivial square non-associative tangle G4
\label{G4}}
\end{center}
\end{figure}%
%EndExpansion

\item[2.] Each hexagon $\Sigma^H_k$, $L^H_k = \Sigma^H_k\cap L \subset
\Sigma^H_k\times [0,1]$ with the bracketing of the 3 ordered sets  
$L^H_k\cap I_j$
looks like figure \ref{G5}, i.e. there is no braiding, association, maxima or
minima for $L^H_k$ and the bracketing respect the decomposition of the bottom of
$\Sigma^H_k$ into two components. The type of such a configuration we will
denote by G5.

\end{description}
See figure \ref{quant1fig1} for an example.

Now we choose a Drinfeld associator as in \cite{Bar-Natan93}. Then we construct
$$V:L(\Sigma) \rightarrow \overline{ch(\Sigma)}$$
 as follows: We associate
to the surface a formal series of chord diagrams as follows: G1-G3 are
mapped as in \cite[section 3]{Bar-Natan93} whereas G4 and G5 are map to
themselves. Since we have no natural way of distinguishing the top and the
bottom of the rectangles, we need to make sure that the associations of G1, G2
and G3 are equivariant with respect to $180^\circ$ rotation. The only real
issue here is the associator. However, by Prop. 3.1. in \cite{Le96a} we know
associators with this symmetry property exits. Alternatively, one can choose
the intervals $I_j$, such that a consistent choice of top and bottom on all
squares and hexagons can be made. In this case we do not need the symmetry
assumption on the associator.

\emph{Claim}\textbf{:} The map $V$ is an
invariant of the link $L$.

This holds since \newline
\begin{enumerate}

\item  Isotopies of links can be assumed to be generic with respect to the
decomposition of the surface into squares and hexagons. So isotopies are
sequences of isotopies in compositions of morphisms, hence leaves $V$ invariant
by the relations R1-R8 of
\cite{Bar-Natan93}, and of moves R10, naturality of the hexagons, see figure
\ref{quant1fig2}, which obviously also leaves
$V$ invariant.

\begin{figure}
\begin{center}
\epsfig{file=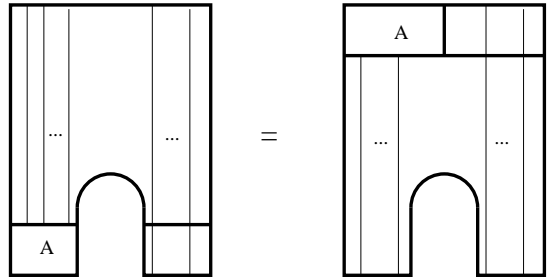,height=7.2cm,width=3.7cm,angle=-90}
\caption{Naturality of the hexagon, R10\label{quant1fig2}}
\end{center}
\end{figure}%

\item  Any two choices of bracketings on the same underlying link give the
same invariant since we can assume by R10 that the bracketings agree on the
hexagons which again reduces the problem to one inside
compositions of morphisms.

\end{enumerate}

Finally, equation (\ref{a230896}) follows from the form of $R,R^{-1}$ as in
the case of quantum group invariants, cf. \cite[theorem 5]{Bar-Natan92}.%
%TCIMACRO{\TeXButton{End Proof}{\endproof}}
%BeginExpansion
$\square$%

%EndExpansion

\begin{corollary}

\label{a230195}If $\Sigma $ is a compact Riemann surface with nonempty
boundary, the space $L^{(n)}(\Sigma )$ is isomorphic to $ch^{(n)}\left(
\Sigma \right) $, the space of link diagrams of order $n$.

\end{corollary}

{\bf Proof:} 
The map $\lambda :ch\left( \Sigma \right) \rightarrow L_{Gr}\left( \Sigma
\right) $ is surjective by definition of $L_{n}$ and injective since $V_{Gr}$
is a left inverse by (\ref{a230896}).

$\square$

\begin{remark}

This contains Goryunov's theorem for the solid torus (see e.g. \cite
{Goryunov95a}) as a special case.

\end{remark}

\begin{remark}

Compare with \cite[theorem 0.1]{t2} where it is shown that for a homotopy $%
\Bbb{R}^{3}$ the theory of Vassiliev invariants is the same as for the usual 
$\Bbb{R}^{3}$.

\end{remark}

\begin{corollary}

\label{a231096}The homomorphism $\lambda :ch\left( \Sigma \right)
\rightarrow L_{Gr}(\Sigma )$ in proposition \ref{b260996} is an isomorphism
of graded Poisson rings.

\end{corollary}

\section{Quantization}

Given a Poisson algebra $\mathcal{A}$ with Poisson bracket $\left\{
.,.\right\} $, a \emph{deformation quantization} of $\mathcal{A}$ is a $\Bbb{%
C}\left[ \left[ h\right] \right] $-linear associative product $*$ on $%
\mathcal{A}\left[ \left[ h\right] \right] $ such that $x*y\text{ mod }h$ is
the given product on $\mathcal{A}$ and $x*y-y*x=h\left\{ x,y\right\} \text{%
mod}h^{2}$ (see \cite{BFFLS}). The product can be written in the form
$x*y=\sum_{n=0}^{\infty }\pi _{n}\left( x,y\right) h^{n}$. We will use the terms
quantization, deformation quantization and star product interchangeably. We
recall that two such star products are equivalent, if there is a
${\Bbb C}[[h]]$-linear automorphism of $\mathcal{A}\left[ \left[ h\right] 
\right]
$, which induces the identity mod $h$ and which is a homomorphism with respect to
the two star products.

In the first subsection we discuss the classical way of quantizing linear
Poisson structures. This applies to the algebra of chord diagrams, but we
show that the resulting star-product does not descend to a star product on
the moduli space of flat $SL_{2}$-connections on a Riemann surface.
Following this discussion we give a geometric construction of a quantization
of the algebra of chord diagrams and verify for $Gl_{n}$ and $SL_{n}$ that
this latter quantization does descend to a deformation quantization of
moduli space.

\subsection{Quantization of chord diagrams as a linear Poisson algebra}

\label{a120996}

\begin{proposition}

Chord diagrams that are connected as a topological space are closed under
the Poisson bracket. The algebra of all chord diagrams is isomorphic to the
polynomial algebra over the set of connected chord diagrams.

\end{proposition}

%TCIMACRO{\TeXButton{Proof}{\proof}}

%BeginExpansion

{\bf Proof:} Obvious from the definitions. $\square$

The proposition says that the algebra of chord diagrams is a linear Poisson
algebra, i.e. a Poisson algebra that is obtained from a Lie algebra $\frak{g}
$ (connected chord diagrams in our case) by extending the Lie bracket as a
derivation to the symmetric algebra $S\frak{g}$ over the given Lie algebra.
Linear Poisson structure can always be quantized in the following way:

Consider the algebra $T\frak{g}\left[ \left[ h\right] \right] /\langle
x\bigotimes y-y\bigotimes x-h\left[ x,y\right] \rangle $. This is an
associative algebra and as a linear space it is isomorphic to $S\frak{g}%
\left[ \left[ h\right] \right] $, it defines a deformation quantization $*$
of $S\frak{g}\left[ \left[ h\right] \right] $.

\begin{lemma}

For $x,y\in \frak{g}$ we have $x*y=xy+\frac{h}{2}\left\{ x,y\right\} $ and $%
x^{2}*y=x^{2}y+h\left\{ x^{2},y\right\} +\frac{h^{2}}{6}\left\{ x,\left\{
x,y\right\} \right\} $.

\end{lemma}

{\bf Proof:} Direct computation. $\square$

Using this one can show that e.g. for $SL_{2}$ this star-product does not
descend to a star product on the moduli space of flat connections:

\begin{example}

On the 4-punctured sphere let $x_{i},i=1,\dots ,4$ denote the generators of
the fundamental group, $\left( 1\right) $ a trivial loop and $\left(
x_{i_{1}}\dots x_{i_{n}}\right) $ the loop given as the corresponding
product of generators. The Cayley-Hamilton theorem implies that the linear
combination of chord diagrams $\left( x_{1}x_{2}x_{1}x_{2}\right) -\left(
x_{1}x_{2}\right) \left( x_{1}x_{2}\right) +\left( 1\right) $ is in the
kernel of the map to moduli space. Taking star product with $\left(
x_{1}x_{3}\right) $ we get 
\begin{multline*}
\left( \left( x_{1}x_{2}x_{1}x_{2}\right) -\left( x_{1}x_{2}\right) \left(
x_{1}x_{2}\right) +\left( 1\right) \right) *\left( x_{1}x_{3}\right) \\
=\left( \left( x_{1}x_{2}x_{1}x_{2}\right) -\left( x_{1}x_{2}\right) \left(
x_{1}x_{2}\right) +\left( 1\right) \right) \left( x_{1}x_{3}\right) \\
+h\left\{ \left( \left( x_{1}x_{2}x_{1}x_{2}\right) -\left(
x_{1}x_{2}\right) \left( x_{1}x_{2}\right) \right) ,\left( x_{1}x_{3}\right)
\right\} \\
+\frac{h^{2}}{6}\left\{ \left( x_{1}x_{2}\right) \left( x_{1}x_{2}\right)
,\left\{ \left( x_{1}x_{2}\right) \left( x_{1}x_{2}\right) ,\left(
x_{1}x_{3}\right) \right\} \right\}
\end{multline*}

The first two lines are zero since the ideal is a Poisson ideal and $\left\{
\left( 1\right) ,\left( x_{1}x_{3}\right) \right\} =0$ but the last
expression does not map to zero on $\mathcal{M}^{SL_{2}}$, hence the star
product does not respect the kernel of the quotient map.

\end{example}

\subsection{Universal quantization of the algebra of chord diagrams}

By universality of the quantization we mean the following: The quantization
descends to the algebra of functions on moduli space in a functorial way
with respect to group homomorphisms that respect the invariant inner product
on the groups. In this subsection we construct the quantization, in the next
subsection we discuss the universality property. The construction uses the
following general facts about filtered algebras:

\begin{proposition}

Let $G=\bigoplus_{n=0}^{\infty }G^{\left( n\right) }$ and $%
H=\bigoplus_{n=0}^{\infty }H^{\left( n\right) }$ be graded spaces with
associated filtered spaces $\overline{G}$ and $\overline{H}$ respectively.
If $\Lambda :\overline{H}\rightarrow \overline{G}$ is linear map respecting
the filtrations such that the associated graded map $\Lambda
_{Gr}:H\rightarrow G$ is an isomorphism, then $\Lambda $ is an isomorphism.
\end{proposition}

Recall from
lemma
\ref{a211096} that for a filtered algebra $F$ such that the associated graded
algebra $F_{Gr}$ is commutative $F_{Gr}$ has a natural structure of a Poisson
algebra. For a graded algebra $G$ we have the projection $x\in G\longmapsto
x^{\left( i\right) }\in G^{\left( i\right) }$, which obviously extends to the
completion $\overline{G}$. Likewise, if we have a $u\in F_{i}$, we denote the
project of $u$ onto $F^{(i)}$ by $u^{(i)}$.

\begin{theorem}

\label{b231096}Let $G$ be a graded Poisson algebra and $F$ a filtered
algebra such that the associated graded algebra is commutative. Assume that $%
\Lambda :F\rightarrow \overline{G}$ is a linear filtered map such that
$\Lambda _{Gr}:F_{Gr}\rightarrow G$ is a Poisson algebra isomorphism. Then
setting 
\[
x_{1}*x_{2}:=h^{-\deg \left( x_{1}x_{2}\right) }\sum_{i=0}^{\infty }\left(
\Lambda \left( \Lambda^{-1} \left( x_{1}\right) .\Lambda^{-1} \left( x_{2}\right)
\right) \right) ^{\left( i\right) }h^{i} 
\]
for $x_i \in G^{(\deg x_i)}$ defines a star product on $G\left[ \left[ h\right]
\right]
$.

\end{theorem}

{\bf Proof:} We have to show that the product is associative and that 
\begin{eqnarray}
x_{1}*x_{2} &=&x_{1}x_{2}\text{ mod }h \\
x_{1}*x_{2}-x_{2}*x_{1} &=&h\left\{ x_{1},x_{2}\right\} \text{ mod }h^{2}
\end{eqnarray}
Associativity holds because the pullback $x_{1}\circ x_{2}:=\Lambda
\left( \Lambda^{-1} \left( x_{1}\right) .\Lambda^{-1} \left( x_{2}\right)
\right) 
$ of an associative multiplication is an associative multiplication and the
coefficients of $h^{p}$ in $\left( x_{1}*x_{2}\right) *x_{3}$ and $%
x_{1}*\left( x_{2}*x_{3}\right) $ are just the coefficients of $\left(
x_{1}\circ x_{2}\right) \circ x_{3}$ and $\left( x_{1}\circ x_{2}\right)
\circ x_{3}$ respectively in $G ^{(\deg (x_{1}x_{2}x_{3})+p)}$.

To see the equations relating the star product to the given product and the
Poisson bracket we compute $\text{ mod }h^2$ 
\begin{eqnarray*}
\lefteqn{x_{1}*x_{2}} \\
&\equiv&(\Lambda(\Lambda^{-1}(x_1).\Lambda^{-1}(x_2)))^{(\deg (x_1x_2))} +
(\Lambda(\Lambda^{-1}(x_1).\Lambda^{-1}(x_2)))^{(\deg (x_1x_2)+1)}h\\
&=&\Lambda_{Gr}(\Lambda_{Gr}^{-1}(x_1)\Lambda_{Gr}^{-1}(x_2)) +
(\Lambda(\Lambda^{-1}(x_1).\Lambda^{-1}(x_2)))^{(\deg (x_1x_2)+1}h\\
&=&x_1x_2 + (\Lambda(\Lambda^{-1}(x_1).\Lambda^{-1}(x_2)))^{(\deg
(x_1x_2)+1)}h \text{ mod }h^2
\end{eqnarray*}
and
\begin{eqnarray*}
\lefteqn{x_{1}*x_{2}-x_2*x_1} \\
&\equiv&(\Lambda(\Lambda^{-1}(x_1).\Lambda^{-1}(x_2)-
              \Lambda^{-1}(x_2).\Lambda^{-1}(x_1)))^{(\deg(x_1x_2)+1)}h\\
&=&\Lambda_{Gr}((\Lambda^{-1}(x_1).\Lambda^{-1}(x_2)-
              \Lambda^{-1}(x_2).\Lambda^{-1}(x_1))^{(\deg(x_1x_2)+1)})h\\
&=&\{x_1,x_2\}h \text{ mod }h^2.
\end{eqnarray*}
We used here that $\Lambda_{Gr}$ is a Poisson homomorphism and the definition
of the Poisson bracket on $F_{Gr}$. $\square$

\begin{remark}

Fixing $\lambda = \Lambda^{-1}_{Gr}$ we see that any two star
products on $G[[h]]$ obtained this way are equivalent.

\end{remark}

Now let us return to the algebra of chord diagrams. Recall from corollary
\ref{a231096} that the map $\lambda =V_{Gr}^{-1}:ch\left( \Sigma \right)
\rightarrow L_{Gr}$ is an isomorphism of graded Poisson rings. Extend
$V:L^{\prime }\rightarrow \overline{ch\left(
\Sigma \right) }$ by continuity to $\overline{V}:\overline{L^{\prime }}%
\rightarrow \overline{ch\left( \Sigma \right) }$, then $\overline{V}$
satisfies the assumptions of theorem \ref{b231096}.

The star product on $ch\left( \Sigma \right) $ obtained from theorem \ref
{b231096} is given as follows: Write $V\left( L\right) =\sum_{D\in ch\left(
\Sigma \right) }\langle D,L\rangle D=\sum_{V\left( L\right) ^{\left(
i\right) }\in ch^{\left( i\right) }\left( \Sigma \right) }V\left(
L\right) ^{\left( i\right) }$ and define a map $F$ from $L^{\prime }\left(
\Sigma \times \left[ 0,1\right] \right) $ to $ch\left( \Sigma \right) \left[
\left[ h\right] \right] $ by 

\[
F\left( L\right) :=\sum_{D\in ch\left( \Sigma \right) }\langle D,L\rangle
Dh^{\deg D}=\sum V\left( L\right) ^{\left( i\right) }h^{i} 
\]

Now the multiplication on $ch\left( \Sigma \right) \left[ \left[ h\right]
\right] $ reads 
\begin{equation}
D_{1}*D_{2}:=h^{-\deg \left( D_{1}D_{2}\right) }F\left( V^{-1}\left(
D_{1}\right) .V^{-1}\left( D_{2}\right) \right)  \label{b120996}
\end{equation}

\begin{theorem}

The product (\ref{b120996}) defines a deformation quantization of $ch\left(
\Sigma \right) $.

\end{theorem}

\begin{remark}

Instead of $\overline{V}$ we could have used $\overline{V_{Gr}}$ and we
would have obtained an isomorphic product.

\end{remark}

\begin{remark}

For the case of closed surfaces computer calculations using Dror Bar-Natan's
mathematica package \cite{Bar-Natan93} suggest that the generalization of
the above construction is nontrivial. Notice that in
the work of Alekseev, Grosse and Schomerus on quantization of moduli space
the case of closed surfaces also presented additional difficulties, in fact
they formulate their result for closed surfaces only as a conjecture \cite
{Alekseev95}.

\end{remark}

\subsection{Quantization of $\mathcal{M}^{G}$}

Here we show that (in contrast to the quantization in section \ref{a120996})
the quantization of the algebra of chord diagrams given above descends to
the moduli space of flat $G$-connections. The proof rests on the following
locality result for the universal Vassiliev invariants constructed above:

\begin{theorem}

Suppose a subspace $W$ in $ch\left( \Sigma \right) $ is defined by local
(''skein'') relations, then it is a (two-sided) ideal with respect to the
star multiplication.

\end{theorem}

{\bf Proof:} 
This follows from the fact that the map $U$ in theorem \ref{b161096} is
multiplicative by definition. Since $U$ is multiplicative, so is $U^{-1}$
(because $U\left( U^{-1}\left( S\right) U^{-1}\left( T\right) \right) =ST$
implies $U^{-1}\left( S\right) U^{-1}\left( T\right) =U^{-1}\left( ST\right)
)$.

Let $D_{1}\in W$ and $D_{2}$ be two chord diagrams, $\Bbb{D}$ some disk on
the surface containing all chords and crossings, $d_{i}$ the chord tangle
obtained by intersecting $D_{i}$ with $\Bbb{D}$, $L=V^{-1}\left(
D_{1}\right) V^{-1}\left( D_{2}\right) $ and $l=L\cap \Bbb{D}$. Recall that $%
V\left( L\right) $ is given by 'closing up' $U\left( l\right) $ on the
surface, analogously $V^{-1}\left( D_{i}\right) $ is given by 'closing up' $%
U^{-1}\left( d_{i}\right) $.

Since the relation $R$ is local we can assume that $d_{1}=\tilde{d}_{1}R\bar{%
d}_{1}$ and $l=\tilde{l}\left( U^{-1}R\otimes \text{id}\right) \bar{l}$
for suitable (chord) tangles $\tilde{d}_{1},\bar{d}_{1},\tilde{l},\bar{l}$.
This implies 
\[
U\left( U^{-1}\left( D_{1}\right) U^{-1}\left( D_{2}\right) \right) =U\left( 
\tilde{l}\right) U\left( U^{-1}R\otimes \text{id}\right) U\left( \bar{l%
}\right) =U\left( \tilde{l}\right) \left( R\otimes \text{id}\right)
U\left( \bar{l}\right) 
\]
which shows that $V\left( V^{-1}\left( D_{1}\right) V^{-1}\left(
D_{2}\right) \right) $ is in $W$. Checking that the powers of $h$ work out
correctly shows that $D_{1}*D_{2}\in W$ so that $W$ is a left ideal. The
same proof shows that it is also a right ideal. 
%TCIMACRO{\TeXButton{End Proof}{\endproof}}
%BeginExpansion
$\square$%

%EndExpansion

Next we show that the star product is compatible with the passage to the
loop algebra (cf. \cite{m1},\cite[proposition 5]{Mattes96a}):

\begin{corollary}

The star product (\ref{b120996}) is compatible with the ideal in $ch\left(
\Sigma \right) $ defined at the end of section 3 of \cite{Mattes96a}.

\end{corollary}

%TCIMACRO{\TeXButton{Proof}{\proof}}

%BeginExpansion

{\bf Proof:} 
The relations defining the ideal are local on the surface. 
%TCIMACRO{\TeXButton{End Proof}{\endproof}}
%BeginExpansion
$\square$%

%EndExpansion

We still have to describe the ideal in the loop algebra which will map to $0$
on moduli space. All functions on $M_{n}^{k}$, invariant under the diagonal $%
GL_{n}$-action, are traces of products. The relations between them are given
by the trace identities for $n\times n$ matrices of \cite[theorem 4.5]
{Procesi76}. The ideal of relations is generated by the fundamental trace
identities which have the form $\sum_{\sigma }a_{\sigma }\text{tr}\mu
_{\sigma }\left( A_{1}\otimes \dots \otimes A_{n}\right) $ where $%
\sigma =\left( i_{1}\dots i_{k}\right) \left( j_{1}\dots j_{h}\right) \dots
\left( t_{1}\dots t_{e}\right) $ is a permutation, the $A_{i}$ are matrices
and 
\[
\mu _{\sigma }\left( A_{1}\otimes \dots \otimes A_{n}\right) =\left(
A_{i_{1}}\dots A_{i_{k}}\right) \otimes \left( A_{j_{1}}\dots
A_{j_{h}}\right) \otimes \dots \otimes \left( A_{t_{1}}\dots
A_{t_{e}}\right) 
\]

In terms of chord diagrams this is again a local relation as in the lemma.

The moduli space $GL_{n}^{k}/GL_{n}$ has as additional invariant functions
only those of the form $\det^{-1}A_{i}$ which can be written as traces of
representations of monomials in the $A_{i}^{-1}$. The only new relations
will be $\det^{-1}A_{i}\det A_{i}=1$. Analogous statements hold for $%
SL_{n}^{k}/SL_{n}$.

This means that all relations defining the kernel of the map $ch\left(
\Sigma \right) \rightarrow \mathcal{F}\left( \mathcal{M}^{G}\right) $ are
local on $\Sigma $. Therefore the above theorem implies

\begin{theorem}

The star product (\ref{b120996}) descends to a star product on $\mathcal{M}%
^{G}$ for $G=GL\left( n,\Bbb{C}\right) ,SL\left( n,\Bbb{C}\right) $.

\end{theorem}

A  discussion of other Lie groups and relations to
other attempts at quantizing the moduli space will appear elsewhere. Here we
only note that we are led to the following

\begin{conjecture}

The star product (\ref{b120996}) descends to a star product on $\mathcal{M}%
^{G}$ for any simple group $G$.

\end{conjecture}

\begin{remark}

Recall that the Poisson homomorphism from the algebra of chord diagrams to
the moduli space of flat connections is universal in the sense that the
following diagram commutes

%TCIMACRO{

%\TeXButton{CD230195b}{\[\begin{CD}

%ch(\Sigma)               @>id>>  ch(\Sigma)\\

%@AAA                                     @AAA\\

%ch(\Sigma)^{G_2} @>\phi ^{*}>>   ch(\Sigma )^{G_1}\\

%@VVfV                                    @VVfV   \\

%\mathcal{F}(\mathcal{M}_{\Sigma}^{G_2})  @>(\phi _*)^{*}>> \mathcal{F}(\mathcal{M}_{\Sigma}^{G_1})

%\end{CD}\]

%} }

%BeginExpansion

\[\begin{CD}
ch(\Sigma)               @>id>>  ch(\Sigma)\\
@AAA                                     @AAA\\
ch(\Sigma)^{G_2} @>\phi ^{*}>>   ch(\Sigma )^{G_1}\\
@VVfV                                    @VVfV   \\
\mathcal{F}(\mathcal{M}_{\Sigma}^{G_2})  @>(\phi _*)^{*}>> \mathcal{F}(\mathcal{M}_{\Sigma}^{G_1})
\end{CD}\]
%
%EndExpansion
for a given group homomorphism $\phi :G_{1}\rightarrow G_{2}$ respecting the
inner products on $G_{1},G_{2}$. The homomorphisms in the diagram are
described in \cite{Mattes96a}. Clearly the quantized version has an
analogous universality property 

%TCIMACRO{

%\TeXButton{CD150996a}{\[\begin{CD}

%ch(\Sigma)[[h]]               @>id>>  ch(\Sigma)[[h]]\\

%@AAA                                     @AAA\\

%ch(\Sigma)^{G_2}[[h]] @>\phi ^{*}>>   ch(\Sigma )^{G_1}[[h]]\\

%@VVfV                                    @VVfV   \\

%\mathcal{F}(\mathcal{M}_{\Sigma}^{G_2})[[h]]  @>(\phi _*)^{*}>> 

%\mathcal{F}(\mathcal{M}_{\Sigma}^{G_1})[[h]]

%\end{CD}\]

%} }

%BeginExpansion

\[\begin{CD}
ch(\Sigma)[[h]]               @>id>>  ch(\Sigma)[[h]]\\
@AAA                                     @AAA\\
ch(\Sigma)^{G_2}[[h]] @>\phi ^{*}>>   ch(\Sigma )^{G_1}[[h]]\\
@VVfV                                    @VVfV   \\
\mathcal{F}(\mathcal{M}_{\Sigma}^{G_2})[[h]]  @>(\phi _*)^{*}>> 
\mathcal{F}(\mathcal{M}_{\Sigma}^{G_1})[[h]]
\end{CD}\]
%EndExpansion
with respect to the star products.

\end{remark}

%\clearpage

\bibliographystyle{alpha}

\bibliography{books,endnote,mattes,phys,previous}

\end{document}